\documentclass[twocolumn,11pt]{article}

\usepackage[utf8]{inputenc}
\usepackage{setspace}
\usepackage{xpatch}
\usepackage{harvard}

\title{ {\large \bf Optimal Control to Limit the Spread of COVID-19 in Italy} }

\author{
\begin{tabular}{c}
Mohamed Abdelaziz Zaitri${^{1, 2}}$, 
Mohand Ouamer Bibi${^{1}}$,          
Delfim F. M. Torres${^{2,*}}$\\      
\end{tabular}\\[0.3cm]
${^{1}}${\it Research Unit LaMOS (Modeling and Optimization of Systems)}\\ 
{\it Department of Operational Research, University of Bejaia, 06000 Bejaia, Algeria}\\[0.3cm]
${^{2}}${\it Center for Research and Development in Mathematics and Applications (CIDMA)}\\ 
{\it Department of Mathematics, University of Aveiro, 3810-193 Aveiro, Portugal}\\[0.3cm]
{\it *Corresponding author: delfim@ua.pt}\\[0.5cm]
(This is a preprint of a paper whose final and definite form\\
is published  by 'Kuwait Journal of Science' (KJS),\\ 
ISSN 2307-4108 (print), ISSN 2307-4116 (online),\\ 
available at \url{https://journalskuwait.org/kjs}.)}

\usepackage[english]{babel}
\usepackage{amsmath}
\usepackage{amssymb}
\usepackage{graphicx}
\usepackage{epstopdf}
\usepackage{url}

\usepackage{fancyhdr}
\usepackage{pgfplots,multicol}
\usepackage{times}
\usepackage[a4paper,margin=25mm]{geometry}
\usepackage[bf,figurename=Fig.,labelsep=period]{caption}
\usepackage{secdot}
\usepackage[compact]{titlesec}
\titleformat{\subsection}
{\normalfont\normalsize}{\thesubsection}{1pt}{~}
\titleformat{\subsubsection}
{\normalfont\normalsize}{\thesubsubsection}{1pt}{~}
\titlespacing{\subsection}{0pt}{0.4cm}{0pt}
\titlespacing{\subsubsection}{0pt}{0.4cm}{0pt}
\titleformat{\section}{\bf}{\thesection}{1pt}{.~}
\newtheorem{remark}{Remark}

\usepackage{caption}
\usepackage{subcaption}

\begin{document}

\date{ }

\twocolumn[
\begin{@twocolumnfalse}
\maketitle
\thispagestyle{empty}
\begin{center}
{\bf Abstract}
\end{center}
\noindent We apply optimal control theory to a generalized SEIR-type model. 
The proposed system has three controls, representing social distancing, 
preventive means, and treatment measures to combat the spread of the COVID-19 pandemic.
We analyze such optimal control problem with respect to real data transmission in Italy. 
Our results show the appropriateness of the model, 
in particular with respect to the number 
of quarantined/hospitalized (confirmed and infected) and recovered individuals. 
Considering the Pontryagin controls, we show how in a perfect world 
one could have drastically diminish the number of susceptible, 
exposed, infected, quarantined/hospitalized, and death individuals, by increasing 
the population of insusceptible/protected.\newline\newline
\vspace{1cm}
\end{@twocolumnfalse}
] 


\section{Introduction}

A severe outbreak of respiratory illness started in Wuhan, 
a city of eleven million people in central China, in December 2019. 
The causative agent was the novel severe acute respiratory 
syndrome coronavirus 2 (SARS-CoV-2), which was identified 
and isolated from a single patient in early January 2020 and subsequently 
verified in sixteen additional patients. The virus is believed to have a zoonotic origin. 
In particular, the Huanan Seafood Market, a live animal and seafood wholesale market in Wuhan, 
was regarded as a primary source of this epidemic, as it is found that 55\% of the first 
four hundred twenty-five confirmed cases were linked to the marketplace. 
Meanwhile, recent comparisons  of the genetic sequences 
of this virus and bat coronaviruses show a 96\% similarity \cite{intro}.

Multiple mathematical models were already presented to predict the dynamics of this pandemic 
at a regional and global level, and some of these models were implemented,
following different methods, to evaluate a strategy for preventive measures:
in \cite{india}, the classical susceptible--infected--recovered (SIR) modeling 
approach \cite{Kermack} was employed to study the parameters of this model for India
while considering different governmental lockdown measures; in \cite{aalpha},
the length of the incubation period of COVID-19 is estimated using confirmed COVID-19 cases 
reported between January 4 and February 24, 2020, from fifty provinces, regions, 
and counties from China; in \cite{bbeta} a model of the outbreak in Wuhan, 
with individual reaction and governmental action (holiday extension, city lockdown, 
hospitalisation and quarantine) is analyzed in the light of the 1918 
influenza pandemic in London; in \cite{telmcen}, 
susceptible--exposed--infectious--recovered (SEIR) modeling is considered
to forecast the COVID-19 outbreak in Algeria by using real data from March 1 
to April 10, 2020; in \cite{Afrique}, a modified SEIR 
model is considered under three intervention scenarios (suppression, mitigation, mildness) 
and simulated to predict and investigate the realities in several African countries: 
South Africa, Egypt, Algeria, Nigeria, Senegal and Kenya. The list of such studies
is long: see, e.g., \cite{LemosP} for a new compartmental epidemiological model for COVID-19 
with a case study of Portugal; \cite{MR4200529} for a fractional (non-integer order) model 
applied to COVID-19 in Galicia, Spain and Portugal; \cite{Zine} for a stochastic 
time-delayed COVID-19 model with application to the Moroccan deconfinement strategy; etc.

In \cite{Peng}, a mathematical system, generalizing the SEIR model,
is presented to analyze the COVID-19 epidemic based on a dynamic 
mechanism that incorporates the intrinsic impact of hidden latent 
and infectious cases on the entire process of the virus transmission. 
The authors of \cite{Peng} validate their model by analyzing data correlation 
on public data of the National Health Commission 
of China from January 20 to February 9, 2020, and 
produce reliable estimates and predictions, revealing 
key parameters of the COVID-19 epidemic. Here, 
we modify the model analyzed in \cite{Peng} in order to consider 
optimal control problems. More precisely, we introduce three control variables 
and combine them with the main parameters of the model of \cite{Peng}. 
Secondly, we analyze a concrete optimal control problem, solving it 
analytically through the celebrated Pontryagin minimum principle \cite{Pontr}.
Moreover, we perform numerical simulations of the spread of COVID-19 
in Italy from September 1 to November 30, 2020.
The model of \cite{Peng} has shown to be a good model 
to describe the reality of China. It's weakness is that it just tries 
to describe a reality but without controlling it. Our main purpose and contribution
here is to include  control measures that allow us to interfere with reality. 
Moreover, we want to illustrate the validity of the model in a different context. 
For this reason, we have considered real data of COVID-19
from Italy instead of China.

The paper is organized as follows. In Section~\ref{sec:2}, we recall 
the generalized SEIR model of \cite{Peng}. Our original results
begin with Section~\ref{sec:3}, where we introduce a generalized 
SEIR control system. An optimal control problem 
is posed and solved analytically in Section~\ref{sec:4}. 
Then, in Section~\ref{sec:5}, we estimate the parameters of the model
using real data of COVID-19 from Italy, and we illustrate the usefulness 
of the proposed optimal control problem through numerical simulations. 
Our results show that the generalized SEIR model of \cite{Peng},
originally considered for China, is also effective with respect to Italy, being
able to model well available real data, while our optimal control approach shows 
clearly the positive and crucial effects of social distancing, preventive means, 
and treatment in the combat of COVID-19. We end with Section~\ref{sec:conc} 
of conclusions.


\section{A generalized SEIR-type model}
\label{sec:2}

The classical SEIR model consists of four compartments: 
susceptible individuals $S(t)$, exposed individuals $E(t)$, infected individuals $I(t)$, 
recovered individuals $R(t)$. This SEIR model is too simplistic to describe COVID-19 
epidemic and new classes need to be included, e.g., Deaths and Quarantined individuals, 
in order to describe the reality. A generalized SEIR-type model for COVID-19 
is proposed by Peng et al. \cite{Peng}, being expressed by a seven-dimensional 
dynamical system as follows:
\begin{equation} 
\label{q1} 
\begin{cases}
\dot{S}(t) = - \dfrac{\beta S(t) I(t)}{N}-\alpha S(t), \\[0.3cm]
\dot{E}(t) = \dfrac{\beta S(t) I(t)}{N} - \gamma E(t), \\[0.3cm]
\dot{I}(t) =  \gamma E(t)-\delta I(t), \\[0.3cm]
\dot{Q}(t) =  \delta I(t)-\lambda(t) Q(t)-\kappa(t) Q(t), \\[0.3cm]
\dot{R}(t) = \lambda(t) Q(t), \\[0.3cm]
\dot{D}(t) = \kappa(t) Q(t), \\[0.3cm]
\dot{P}(t) = \alpha S(t), 
\end{cases}
\end{equation}
subject to fixed initial conditions 
\begin{equation}
\label{eq:ic}
\begin{gathered}
S(0) = S_{0}, \ E(0)=E_0, \ I(0) = I_{0},\ Q(0) = Q_{0},\\ 
R(0)=R_0,\ D(0) = D_{0}, \ P(0) = P_{0}. 
\end{gathered}
\end{equation}
Here, the population is divided into susceptible individuals $S(t)$, 
exposed individuals $E(t)$, infected individuals $I(t)$, 
quarantined/hospitalized individuals (confirmed and infected) $Q(t)$,  
recovered individuals $R(t)$, death individuals $D(t)$, 
and insusceptible individuals (protected population) $P(t)$. 
It follows from \eqref{q1} that 
$$
\dot{S}(t) + \dot{E}(t) + \dot{I}(t) + \dot{Q}(t) + \dot{R}(t) + \dot{D}(t) + \dot{P}(t) = 0,
$$
so that  
$$
S(t) + E(t) + I(t) + Q(t) + R(t) + D(t) + P(t)
$$ 
is constant along time $t$. This constant will be denoted by $N$, 
being determined by the initial conditions \eqref{eq:ic}: 
$$
N := S_{0}+E_0+I_{0}+Q_{0}+R_0+D_{0}+P_{0}.
$$
The constant parameters $\alpha$, $\beta$, $\gamma$ and $\delta$ 
represent, respectively, the protection rate, infection rate, 
inverse of the average latent time, 
and the rate at which infectious people enter in quarantine,
and they have the dimension of time$^{-1}$ (day$^{-1}$).
The recovery and mortality rates, respectively $\lambda$ and $\kappa$, 
are time-dependent analytical functions defined by
\begin{equation} 
\label{eq:lambda}
\lambda(t):=\frac{\lambda_1}{1+e^{-\lambda_2(t-\lambda_3)}}
\end{equation}
and
\begin{equation} 
\label{eq:kappa}
\kappa(t):=\frac{\kappa_1}{e^{\kappa_2(t-\kappa_3)}+e^{-\kappa_2(t-\kappa_3)}},
\end{equation}
where the parameters $\lambda_1$, $\lambda_2$, $\lambda_3$,
$\kappa_1$, $\kappa_2$ and $\kappa_3$ are determined empirically
from real data. Note that $\lambda_1$, $\lambda_2$, $\kappa_1$ and $\kappa_2$ 
have the dimension of time$^{-1}$ (day$^{-1}$), 
while $\lambda_3$ and $\kappa_3$ have the dimension of time (day).

\begin{remark}
The basic reproduction number is usually computed for autonomous systems, 
when the right-hand side of the system does not depend explicitly 
on time $t$ \cite{MR4185283,MR1950747}.
Here, system \eqref{q1} depends on \eqref{eq:lambda} 
and, therefore, it is a non-autonomous system. In this case, 
we are not aware of a valid method to compute the basic reproduction number. 
\end{remark}


\section{Formulation of the Problem}
\label{sec:3}

\begin{figure}
\begin{center}
\includegraphics[scale=0.15]{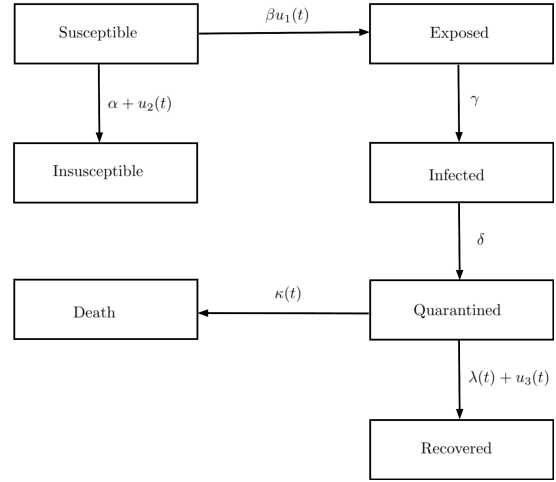}
\caption{Schematic diagram of the generalized SEIR-type control system \eqref{contsys}.}
\label{fig:diagram}
\end{center}
\end{figure}

We introduce three time-dependent controls 
to model \eqref{q1} of \cite{Peng}:
\begin{enumerate}
\item[$-$] control $u_1(t)$, representing the effect of social distancing;
\item[$-$] control $u_2(t)$, representing the effect of preventive means;
\item[$-$] control $u_3(t)$, representing the effect of treatment.
\end{enumerate}
Mathematically, we have the control system
\begin{equation} 
\label{contsys} 
\begin{cases}
\dot{S}(t) = - \dfrac{\beta u_1 (t)S(t) I(t)}{N}-(\alpha+u_2(t)) S(t),\\[0.3cm]
\dot{E}(t) = \dfrac{\beta u_1(t) S(t) I(t)}{N} - \gamma E(t),\\[0.3cm]
\dot{I}(t) =  \gamma E(t)-\delta I(t), \\[0.3cm]
\dot{Q}(t) =  \delta I(t)-(\lambda(t) + u_3(t)) Q(t)-\kappa(t) Q(t), \\[0.3cm]
\dot{R}(t) = (\lambda(t) + u_3(t)) Q(t), \\[0.3cm]
\dot{D}(t) = \kappa(t) Q(t),\\[0.3cm]
\dot{P}(t) = (\alpha+u_2(t)) S(t),
\end{cases}
\end{equation}
subject to initial conditions \eqref{eq:ic}. 
We call \eqref{contsys} the generalized SEIR-type control model.
A schematic diagram of our control system is given in Figure~\ref{fig:diagram}.

\section{Optimal Control}
\label{sec:4}

We consider the generalized SEIR control model \eqref{contsys} 
and formulate an optimal control problem to determine the strategy 
$u(t)=(u_1(t),u_2(t),u_3(t))$, over a fixed interval
of time $[0, t_f ]$, that minimizes the cost functional
\begin{multline}
\label{ju}
J(u)=\int\limits_{0}^{t_f} 
\left(w_1 \frac{\beta u_1 S(t)I(t)}{N}-w_2R(t)\right.\\
\left.-w_3P(t) +v_1 \frac{u_1^2}{2}+v_2\frac{u_2^2}{2}
+v_3\frac{u_3^2}{2}\right) dt,
\end{multline}
where  $t_f$ represents the final time of the period under study 
and the constants $w_1$, $w_2$, $w_3$, $v_1$, $v_2$ and $v_3$ 
represent the weights associated with the total number of new infections, 
the number of recovered individuals, the number of insusceptible individuals,  
and the costs associated with the controls $u_1$, $u_2$ and $u_3$, respectively. 
The controls $u$ are Lebesgue measurable and bounded: 
\begin{multline}
u(t)\in \Gamma:=\left\{\mu=(\mu_1,\mu_2,\mu_3) 
\in \mathbb{R}^3 :\right.\\ 
\left. u_{i\min}\leq \mu_i\leq u_{i\max},\ i=1,2,3\right\}.
\end{multline}
The intervals $[u_{i\min}, u_{i\max}]$ also translate the fact 
that there are limitations to the effects of social distancing, 
the preventive means and the treatment rate. Let 
\begin{equation*}
\begin{split}
x(t)&=(x_1(t),\ldots,x_7(t))\\
&=(S(t),E(t),I(t),Q(t),R(t),D(t),P(t))\\
&\in\mathbb{R}^7. 
\end{split}
\end{equation*}
The optimal control problem consists to find the optimal trajectory $\tilde{x}$ associated 
with the optimal control $\tilde{u}\in L^1$, $\tilde{u}(t) \in \Gamma$, 
satisfying the control system \eqref{contsys}, the initial conditions
\begin{equation}
\label{si}
x(0)=(S_0,E_0,I_0,Q_0,R_0,D_0,P_0)
\end{equation}
and giving minimum value to \eqref{ju}.

The existence of an optimal control $\tilde{u}$ and associated optimal 
trajectory $\tilde{x}$ comes from the convexity of the integrand of the cost 
functional \eqref{ju} with respect to control $u$ and the Lipschitz property 
of the state system with respect to state variables $x$ (see  \cite{existence} 
for existence results of optimal solutions). 
According to the Pontryagin Minimum Principle \cite{Pontr}, 
if $\tilde{u}\in L^1$ is optimal for problem \eqref{contsys}--\eqref{si} 
and fixed final time $t_f$, then there exists $\psi\in AC([0,t_f];\mathbb{R}^7)$,  
$\psi(t) = (\psi_1(t), \ldots , \psi_7(t))$, called the adjoint vector, 
such that
\begin{equation*}
\begin{cases}
\dot{x}=\displaystyle \frac{\partial H}{\partial \psi},\\[0.3cm]
\dot{\psi}=-\displaystyle \frac{\partial H}{\partial x},
\end{cases}
\end{equation*}
where the Hamiltonian $H$ is defined by
\begin{multline*}
H(x,u,\psi )=\frac{w_1u_1\beta x_1 x_3}{N}\\
-w_2 x_5-w_3 x_7 +\sum\limits_{i=1}^3   v_i\frac{u_i^2}{2}\\
+\psi^T   \left( Ax+\left(\sum\limits_{i=1}^2
b_i\Lambda_i x \Phi_i+f(x)^T \Phi_3 \right)u  \right)
\end{multline*}
with
\begin{equation*}
\begin{split}
f(x)&=(f_1(x) \ f_2(x)\  \ 0\ \ 0 \ \ 0 \ \ 0 \ \ 0),\\
f_1(x)&=\frac{-\beta x_1x_3}{N},\\\
f_2(x)&=\frac{\beta x_1x_3}{N},\\
b_1&=(-1\  0\  0 \ 0\  0 \ 0 \ 0)^T,\\ 
b_2&=(0\  0\  \ 0 \ -1\  1 \ 0 \ 0)^T,\\
\Lambda_1&=(1\ 0 \ 0 \ 0\ 0 \ 0 \ 0),\\\
\Lambda_2&=(0\ 0 \ 0 \ 1\ 0 \ 0 \ 0),\\
\Phi_1&=(0\ 1\ 0),\\ 
\Phi_2&=(0\ 0\ 1),\\
\Phi_3&=(1\ 0\ 0),
\end{split}
\end{equation*}
$$
A=\left( 
{\begin{array}{ccccccc}
-\alpha & 0 & 0 &0& 0 & 0 &0\\
0 & -\gamma & 0 & 0& 0 & 0 &0\\
0 & \gamma & -\delta & 0& 0 & 0 &0\\
0 & 0 & \delta & -\lambda(t)-\kappa(t)& 0 & 0 &0\\
0 & 0 & 0 & \lambda(t)& 0 & 0 &0\\
0 & 0 & 0 & \kappa(t)& 0 & 0 &0\\
\alpha & 0 & 0 &0& 0 & 0 &0
\end{array}} 
\right).
$$
The minimality condition
\begin{equation}
\label{ghi}
H(\tilde{x}(t),\tilde{u}(t),\tilde{\psi}(t))
=\min\limits_{u\in\Gamma}H(\tilde{x}(t),u,\tilde{\psi}(t))
\end{equation}
holds almost everywhere on $[0, t_f ]$. Moreover, the 
transversality conditions 
$$
\tilde{\psi}_i(t_f)=0,\quad  i=1,\ldots,7,
$$ 
hold. Solving the minimality condition \eqref{ghi} 
on the interior of the set of admissible controls $\Gamma$ gives
\begin{gather*}
\tilde{u}(t)=\left(\frac{\beta \tilde{x}_1(t)   
\tilde{x}_3(t)\left(\tilde{\psi}_1(t)-\tilde{\psi}_2(t)
-w_1\right)}{Nv_1},\right.\\
\frac{\tilde{x}_1(t)  \left(\tilde{\psi}_1(t)
-\tilde{\psi}_7(t)\right)}{v_2}, \\
\left.\frac{\tilde{x}_4(t)
\left(\tilde{\psi}_4(t)-\tilde{\psi}_5(t)\right)}{v_3}\right),
\end{gather*}
where the adjoint functions satisfy
\begin{equation}
\label{eq:adj:cond}
\begin{cases}
\dot{\tilde{\psi}}_1  =-\displaystyle 
\frac{\tilde{u}_1\beta  \tilde{x}_3}{N^2} 
\left( \tilde{x}_2
+ \tilde{x}_3+ \tilde{x}_4+\tilde{x}_5+ \tilde{x}_6+ \tilde{x}_7\right)\\
\qquad \times \left(w_1 - \tilde{\psi}_1+ \tilde{\psi}_2\right)
+(\alpha+\tilde{u}_2) (\tilde{\psi}_1 - \tilde{\psi}_7) ,\\  
\dot{\tilde{\psi}}_2=\displaystyle \frac{\tilde{u}_1\beta   \tilde{x}_1
\tilde{x}_3\left(w_1-\tilde{\psi}_1+\tilde{\psi}_2\right)}{N^2}
+\gamma (\tilde{\psi} _2-\tilde{\psi}_ 3),\\   
\dot{\tilde{\psi}}_3=
-\displaystyle \frac{\tilde{u}_1\beta \tilde{x}_1}{N^2}
\left(\tilde{x}_2 + \tilde{x}_3+ \tilde{x}_4+\tilde{x}_5+ \tilde{x}_6+ \tilde{x}_7\right)\\
\qquad \times \left(w_1-\tilde{\psi}_1+\tilde{\psi}_2\right)
+\delta (\tilde{\psi}_3-\tilde{\psi}_4),\\      
\dot{\tilde{\psi}}_4
=\displaystyle \frac{\tilde{u}_1\beta\tilde{x}_1\tilde{x}_3
\left(w_1-\tilde{\psi}_1+\tilde{\psi}_2\right)}{N^2}\\
\qquad + \kappa(t) (\tilde{\psi}_4-\tilde{\psi}_6)
+\left(\lambda(t)+\tilde{u}_3\right)(\tilde{\psi}_4-\tilde{\psi}_5),\\
\dot{\tilde{\psi}}_5=\displaystyle \frac{ \tilde{u}_1\beta  \tilde{x}_1  
\tilde{x}_3 (w_1-\tilde{\psi}_1+\tilde{\psi}_2)}{   N^2}+w_2,\\
\dot{\tilde{\psi}}_6=\displaystyle \frac{ \tilde{u}_1\beta  \tilde{x}_1  
\tilde{x}_3 (w_1-\tilde{\psi}_1+\tilde{\psi}_2)}{   N^2},\\
\dot{\tilde{\psi}}_7=\displaystyle \frac{ \tilde{u}_1\beta  \tilde{x}_1  
\tilde{x}_3 (w_1-\tilde{\psi}_1+\tilde{\psi}_2)}{   N^2}+w_3.
\end{cases}
\end{equation}

Note that we have obtained 
an analytical explicit expression for the controls 
$\tilde{u}_1(t), \tilde{u}_2(t)$ and $\tilde{u}_3(t)$,
\begin{equation}
\label{contr123}
\begin{gathered}
\tilde{u}_1(t)=\frac{\beta \tilde{x}_1(t)   
\tilde{x}_3(t)\left(\tilde{\psi}_1(t)-\tilde{\psi}_2(t)
-w_1\right)}{Nv_1},\\
\tilde{u}_2(t)=\frac{\tilde{x}_1(t)  \left(\tilde{\psi}_1(t)
-\tilde{\psi}_7(t)\right)}{v_2}, \\
\tilde{u}_3(t)=\frac{\tilde{x}_4(t)
\left(\tilde{\psi}_4(t)-\tilde{\psi}_5(t)\right)}{v_3},
\end{gathered}
\end{equation}
but we do not have the controls in open-loop (because they depend 
on the state variables $\tilde{x}$ and adjoint variables $\tilde{\psi}$). 
To plot $\tilde{u}(t)$ as a function of $t$ we need to solve numerically 
system \eqref{contsys} and \eqref{eq:adj:cond} to know the expressions 
for $\tilde{x}$ and $\tilde{\psi}$ and be able to obtain the controls 
$u_i$, $i=1,2,3$, in agreement with \eqref{contr123}. This is done numerically
in next section. For more on numerical approaches to solve optimal control
problems, we refer the reader to \cite{MR3758014,MR3953217} and references therein.


\section{Numerical Results}
\label{sec:5}

Now, our aim is to find optimal controls to limit 
the spread of the epidemic of COVID-19 in Italy, by reducing 
the number of new infections and by increasing insusceptible individuals
and the percentage of those recovered, while reducing the cost during 
the period of three months starting from September 1, 2020.  
All numerical computations were performed in the 
numeric computing environment \textsf{MATLAB R2019b} 
using the medium order method and numerical 
interpolation \cite{MR1433374}. 
The rest of the preliminary conditions and real data were taken 
and computed from the database
\url{https://raw.githubusercontent.com/pcm-dpc/COVID-19/master/dati-regioni/dpc-covid19-ita-regioni.csv}.
The real data for COVID-19 pandemic in Italy, 
for September and October 2020, is summarized
in appendix: see Appendix~A for recovered individuals, 
Appendix~B for deaths, and Appendix~C 
for quarantined individuals.

\begin{figure}
\begin{center}
\includegraphics[scale=0.38]{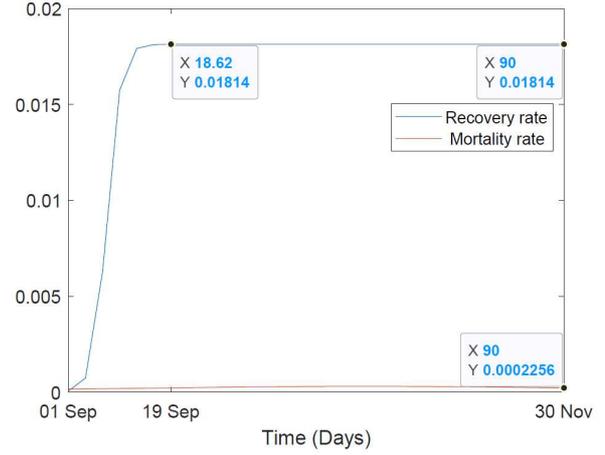}
\caption{The recovery and mortality rates \eqref{eq:lambda}
and \eqref{eq:kappa} for the case of Italy (Section~\ref{sec:5}).}
\label{rmm}
\end{center}
\end{figure}

The parameters $\alpha$, $\beta$, $\gamma$, $\delta$,
$(\kappa_1,\kappa_2,\kappa_3)$ and $(\lambda_1,\lambda_2,\lambda_3)$ 
were fitted in the least square sense. In Figure~\ref{rmm}, we plot 
functions $\lambda(t)$ \eqref{eq:lambda} and $\kappa(t)$ \eqref{eq:kappa} 
by considering the initial guess $\alpha =0.06$, $\beta=1$, $\gamma=5$, 
$\delta=0.5$, $(\lambda_1,\lambda_2,\lambda_3)=(0.01, 0.1, 10)$ 
and $(\kappa_1,\kappa_2,\kappa_3)=(0.001,0.001,10)$, respectively. 

The parameters of the generalized SEIR model \eqref{q1} 
were computed simultaneously by the nonlinear 
least-squares solver \cite{Cheynet}. These parameters, 
during the period under study, were found as follows: 
$\alpha =1.1775\times 10^{-7}$,  
$\beta=3.97$, $\gamma=0.0048$, $\delta=0.1432$, 
$(\lambda_1,\lambda_2,\lambda_3)=(0.0181, 0.8111, 6.9882)$ 
and $(\kappa_1,\kappa_2,\kappa_3)=(0.00062,0.0233,54.0351)$.
For the optimal control problem of Section~\ref{sec:4},
we further fixed $w_i =v_i =1$, $u_{1\min}=0.1$, $u_{j\min}=0$, 
$u_{i\max}=1$, $i=1,2,3$, $j=1,2$. 

\begin{figure*}
\centering
\begin{subfigure}[b]{0.49\textwidth}
\centering
\includegraphics[scale=0.57]{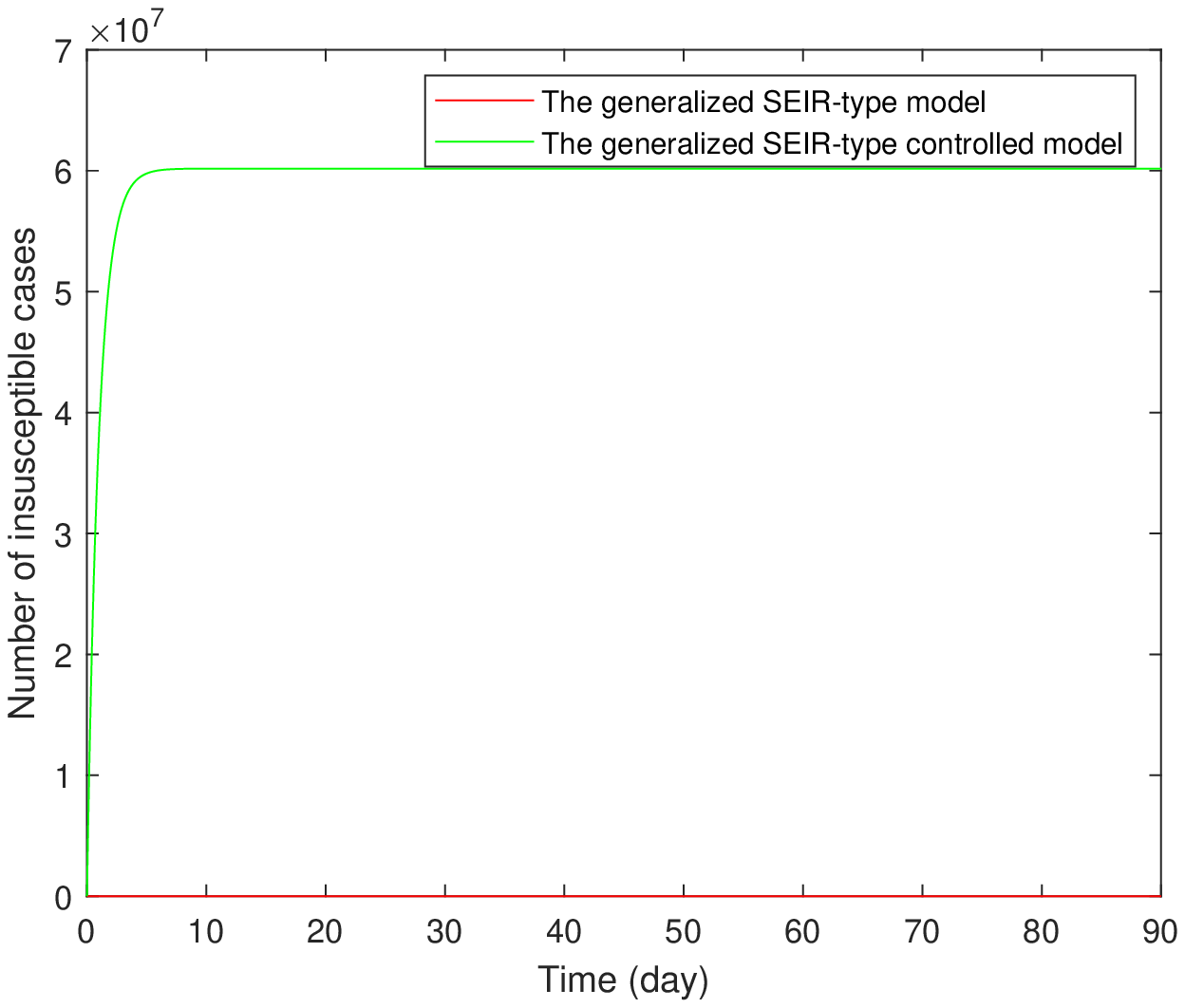}
\caption{$P(t)$}
\label{fig1:a}
\end{subfigure}
\hfill
\begin{subfigure}[b]{0.49\textwidth}
\centering
\includegraphics[scale=0.57]{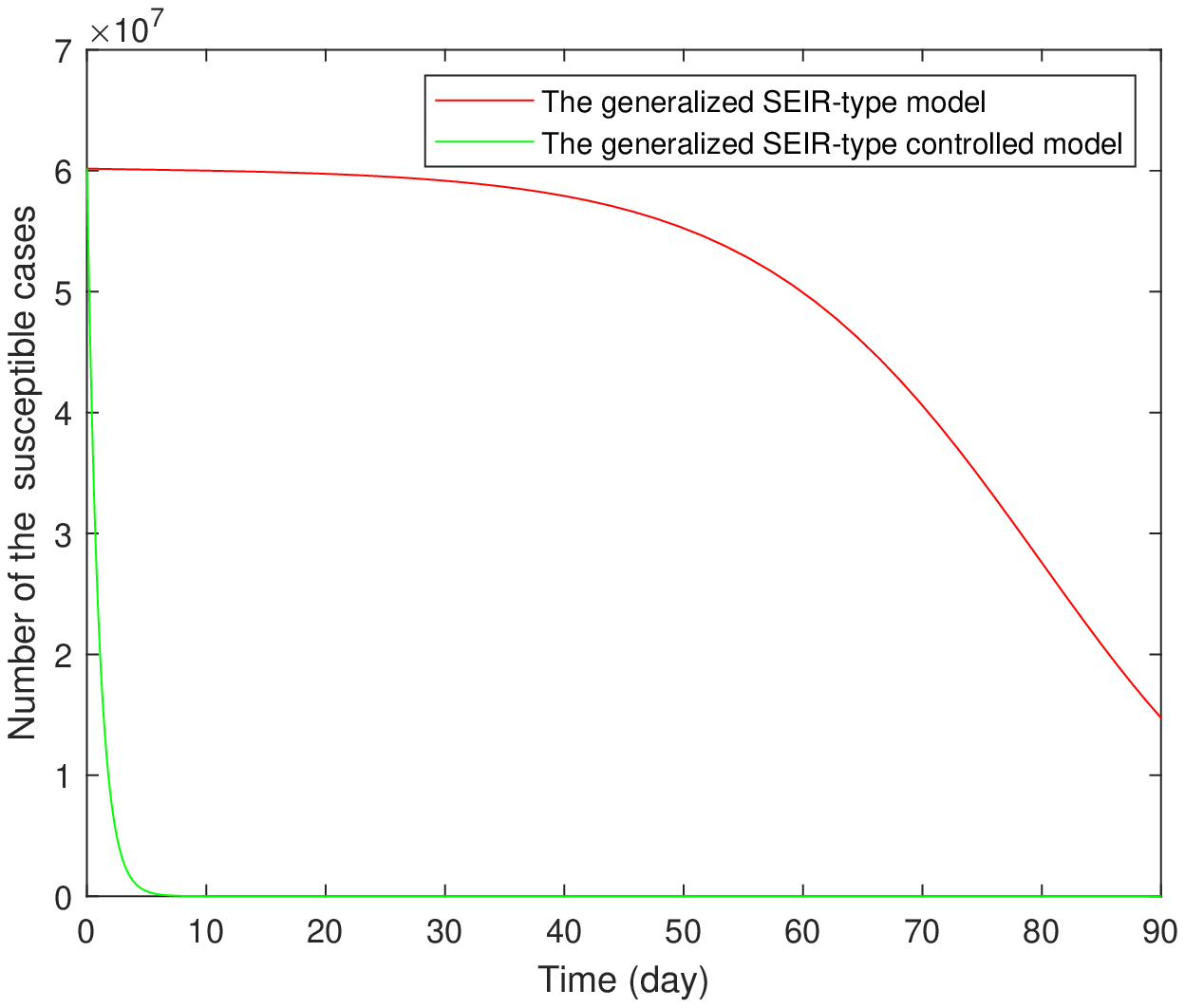}
\caption{$S(t)$}
\label{fig1:b}
\end{subfigure}
\hfill
\begin{subfigure}[b]{0.49\textwidth}
\centering
\includegraphics[scale=0.57]{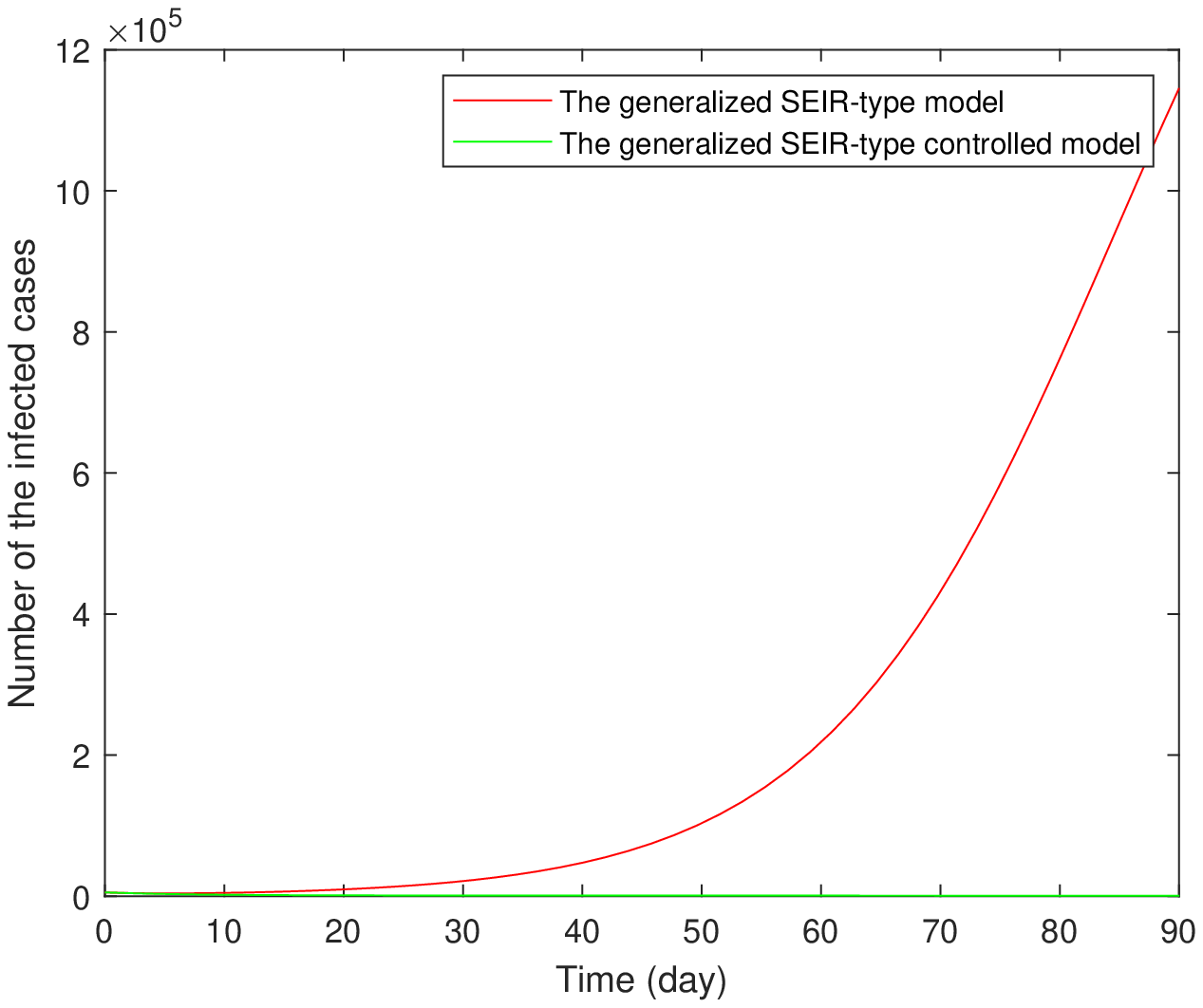}
\caption{$I(t)$}
\label{fig1:c}
\end{subfigure}
\hfill
\begin{subfigure}[b]{0.49\textwidth}
\centering
\includegraphics[scale=0.57]{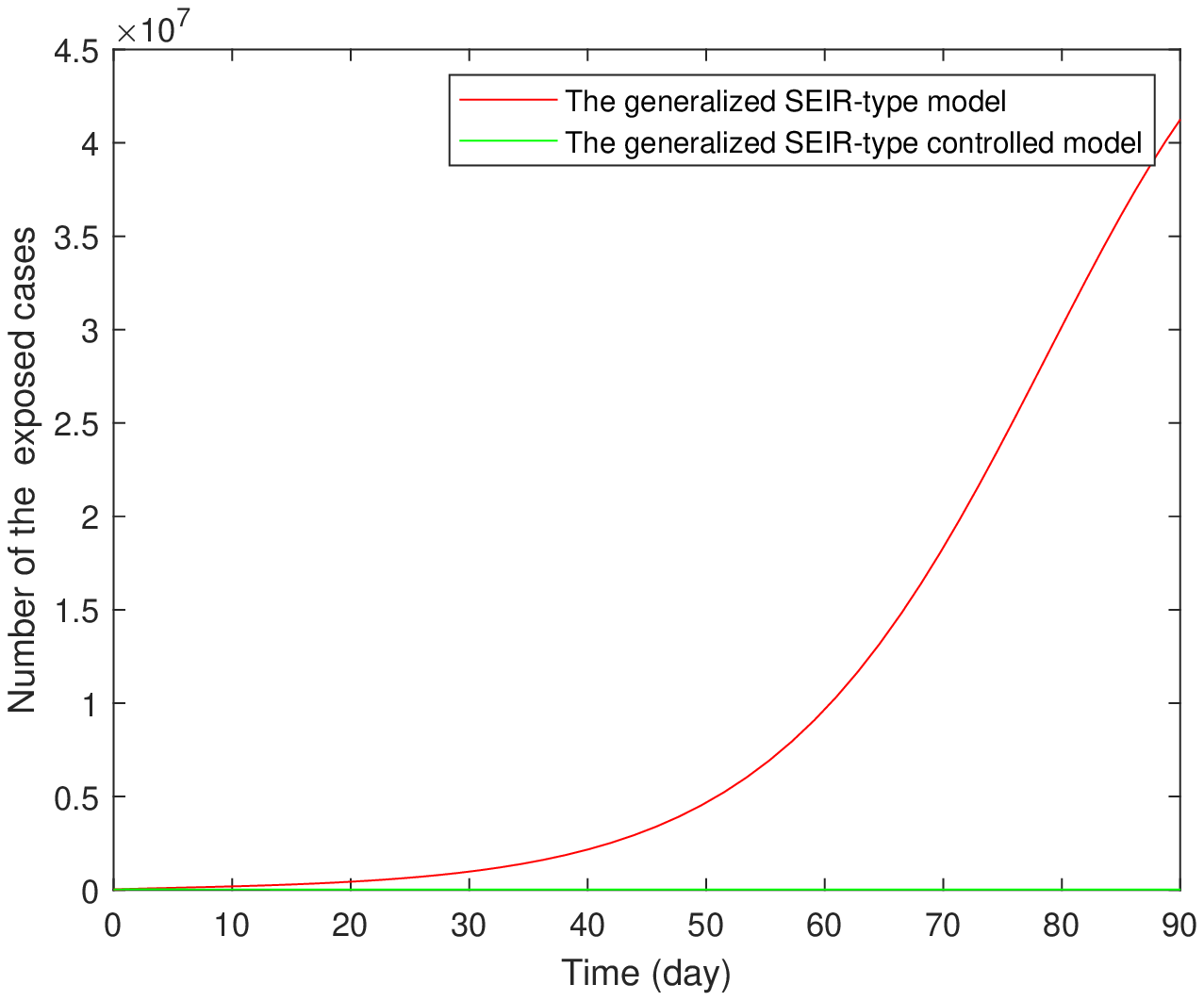}
\caption{$E(t)$}
\label{fig1:d}
\end{subfigure}
\hfill
\begin{subfigure}[b]{0.49\textwidth}
\centering
\includegraphics[scale=0.57]{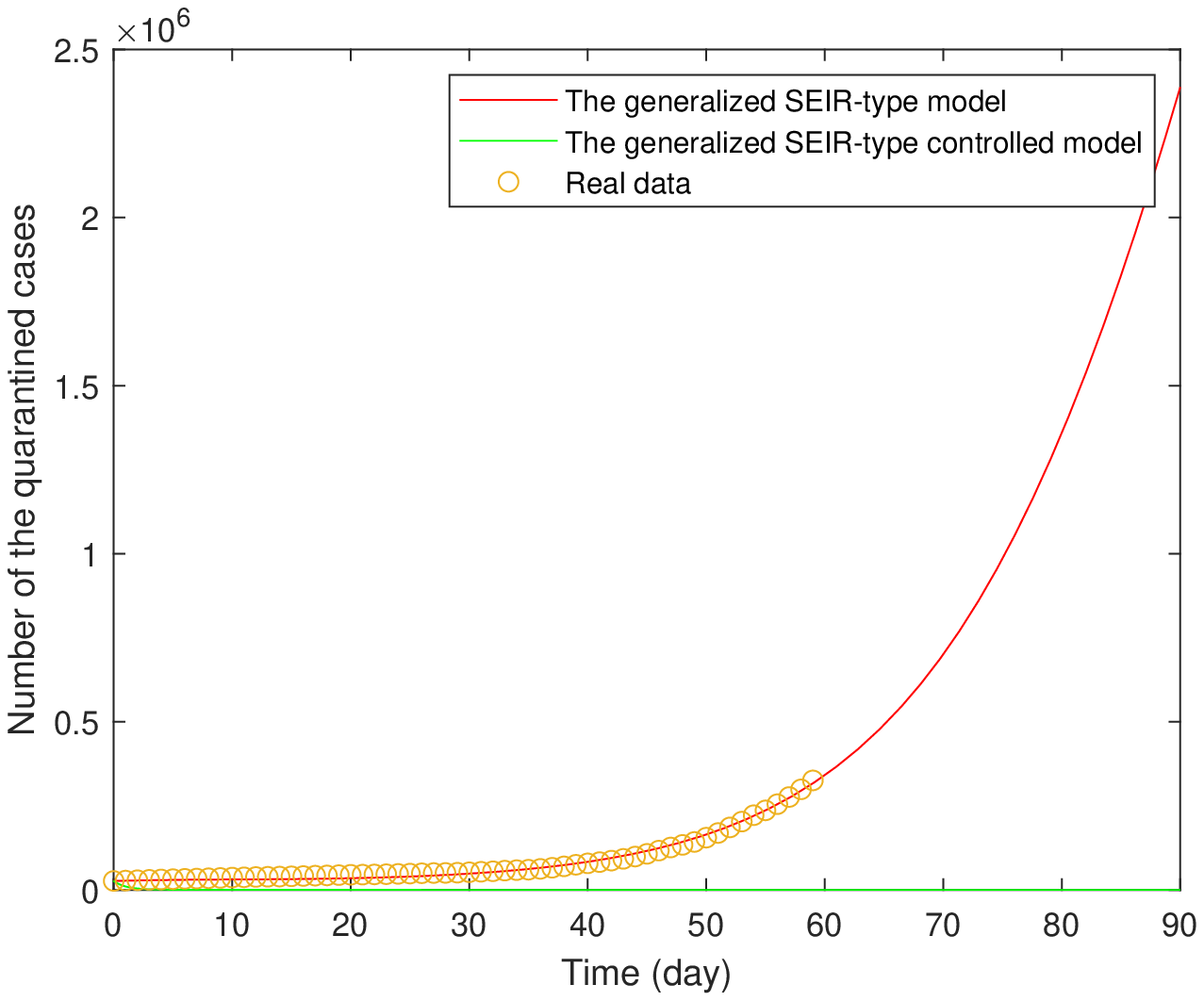}
\caption{$Q(t)$}
\label{fig1:e}
\end{subfigure}
\hfill
\begin{subfigure}[b]{0.49\textwidth}
\centering
\includegraphics[scale=0.57]{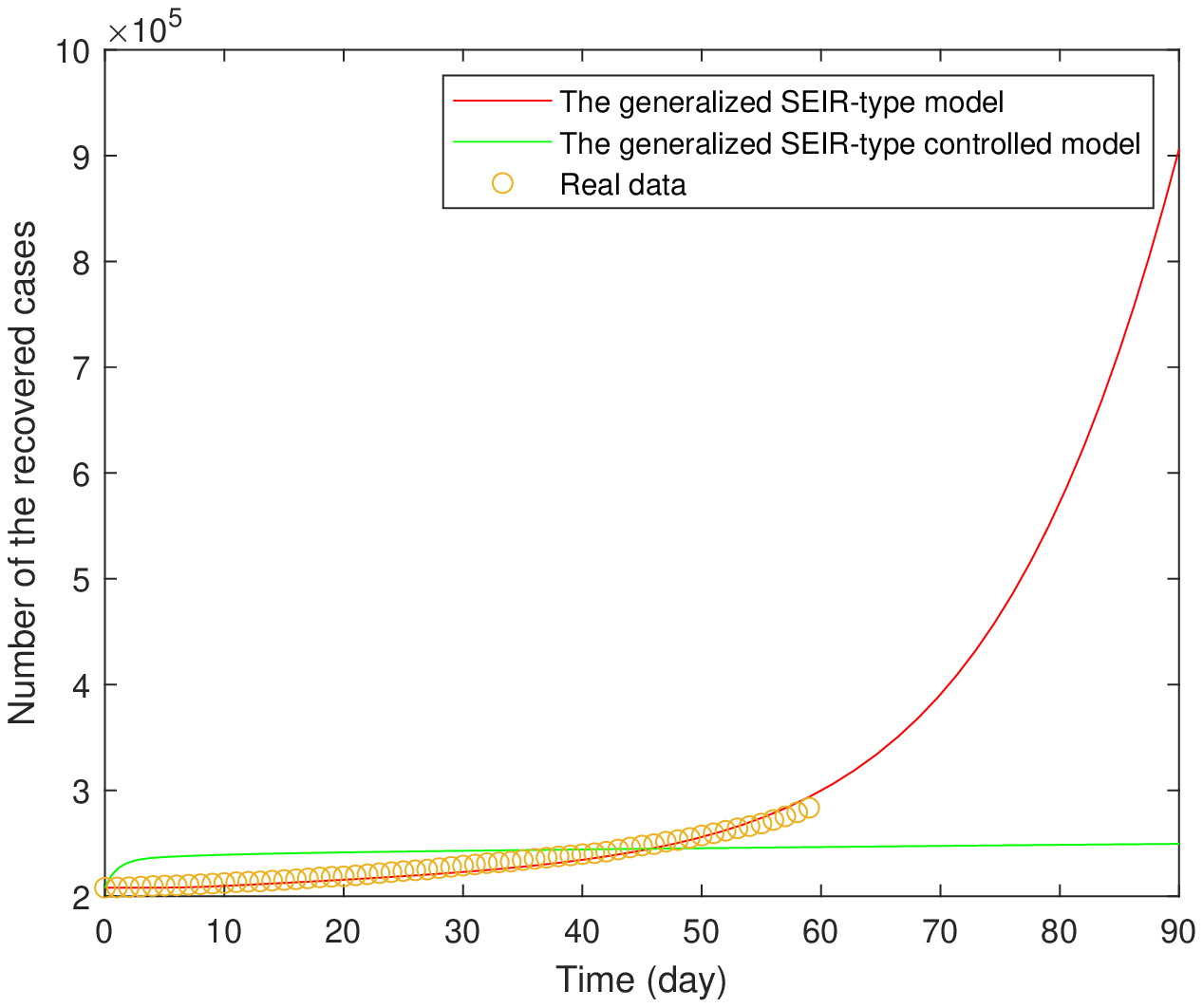}
\caption{$R(t)$}
\label{fig1:f}
\end{subfigure}
\caption{Predictions for Italy from the generalized SEIR model \eqref{q1},
in red, the generalized SEIR control system 
\eqref{contsys} under optimal controls, in green,
between Sept. 1 and Nov. 30, 2020, versus
available real data of quarantined and recovered from Sept. 1 
to Oct. 31, 2020, in orange.}
\label{fig1}
\end{figure*}

\begin{figure}
\begin{center}
\includegraphics[scale=0.57]{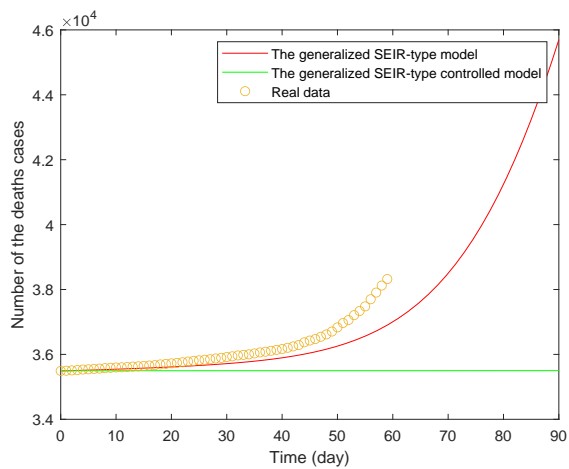}
\caption{Predictions for Italy 
from the generalized SEIR model \eqref{q1},
in red, the generalized SEIR control system 
\eqref{contsys} under optimal controls, in green,
between Sept. 1 and Nov. 30, 2020, versus
available real data of deaths 
from Sept. 1 to Oct. 31, 2020, in orange.}
\label{fig1:g}
\end{center}
\end{figure}

\begin{figure*}
\centering
\begin{subfigure}[b]{0.49\textwidth}
\centering
\includegraphics[scale=0.57]{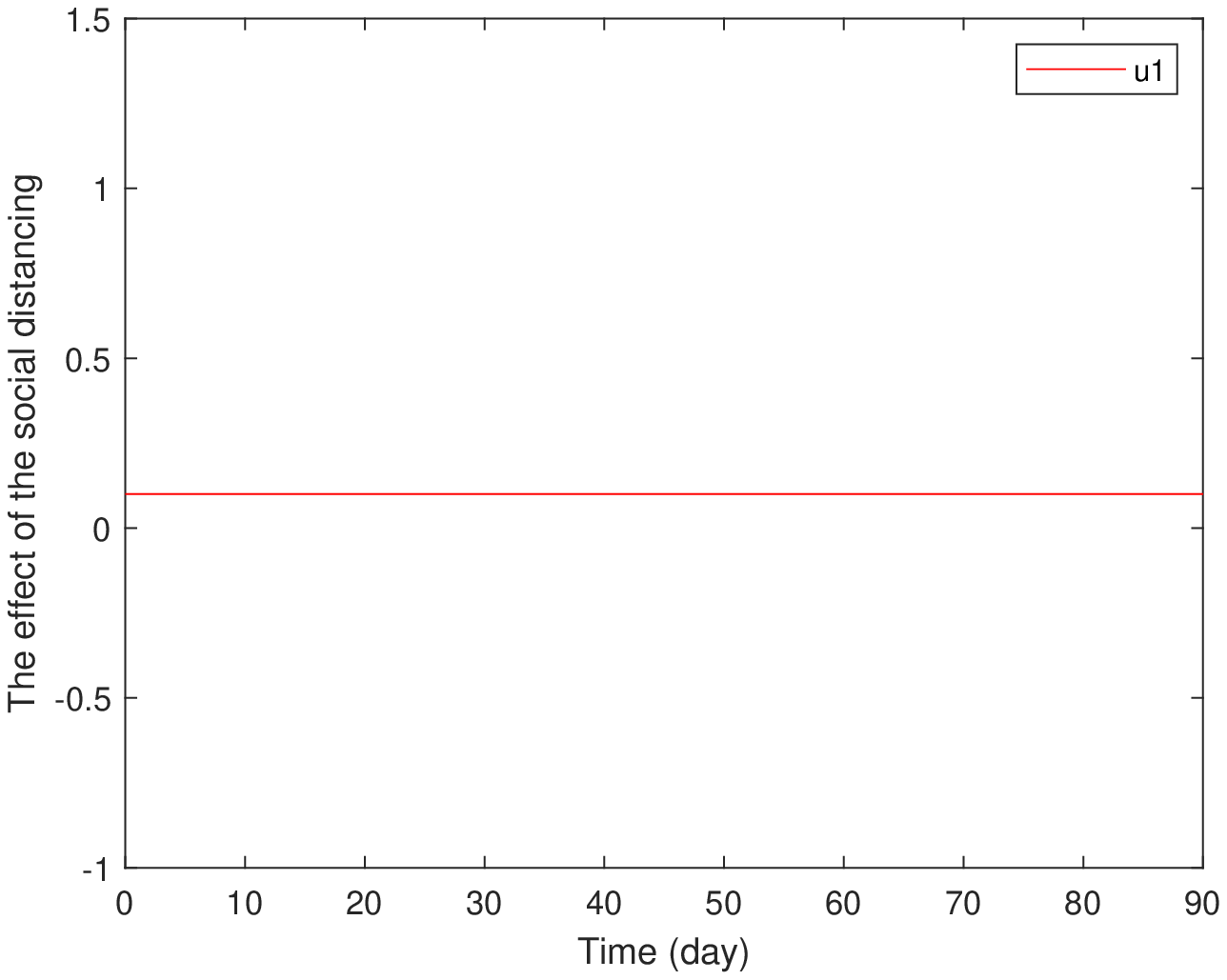}
\caption{$u_1(t)$ (social distancing)}
\label{fig2:c1}
\end{subfigure}
\hfill
\begin{subfigure}[b]{0.49\textwidth}
\centering
\includegraphics[scale=0.57]{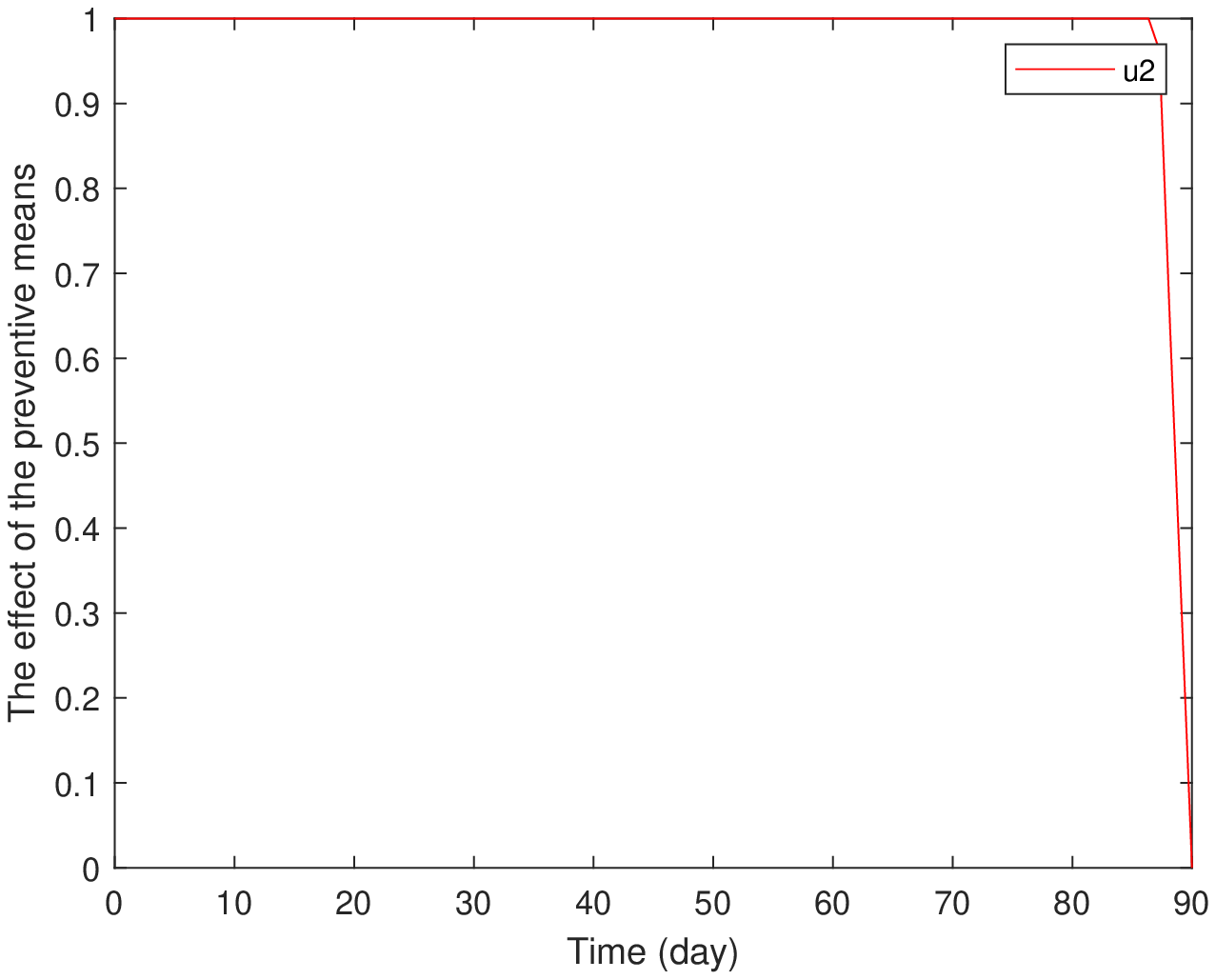}
\caption{$u_2(t)$ (preventive means)}
\label{fig2:c2}
\end{subfigure}
\hfill
\begin{subfigure}[b]{0.49\textwidth}
\centering
\includegraphics[scale=0.57]{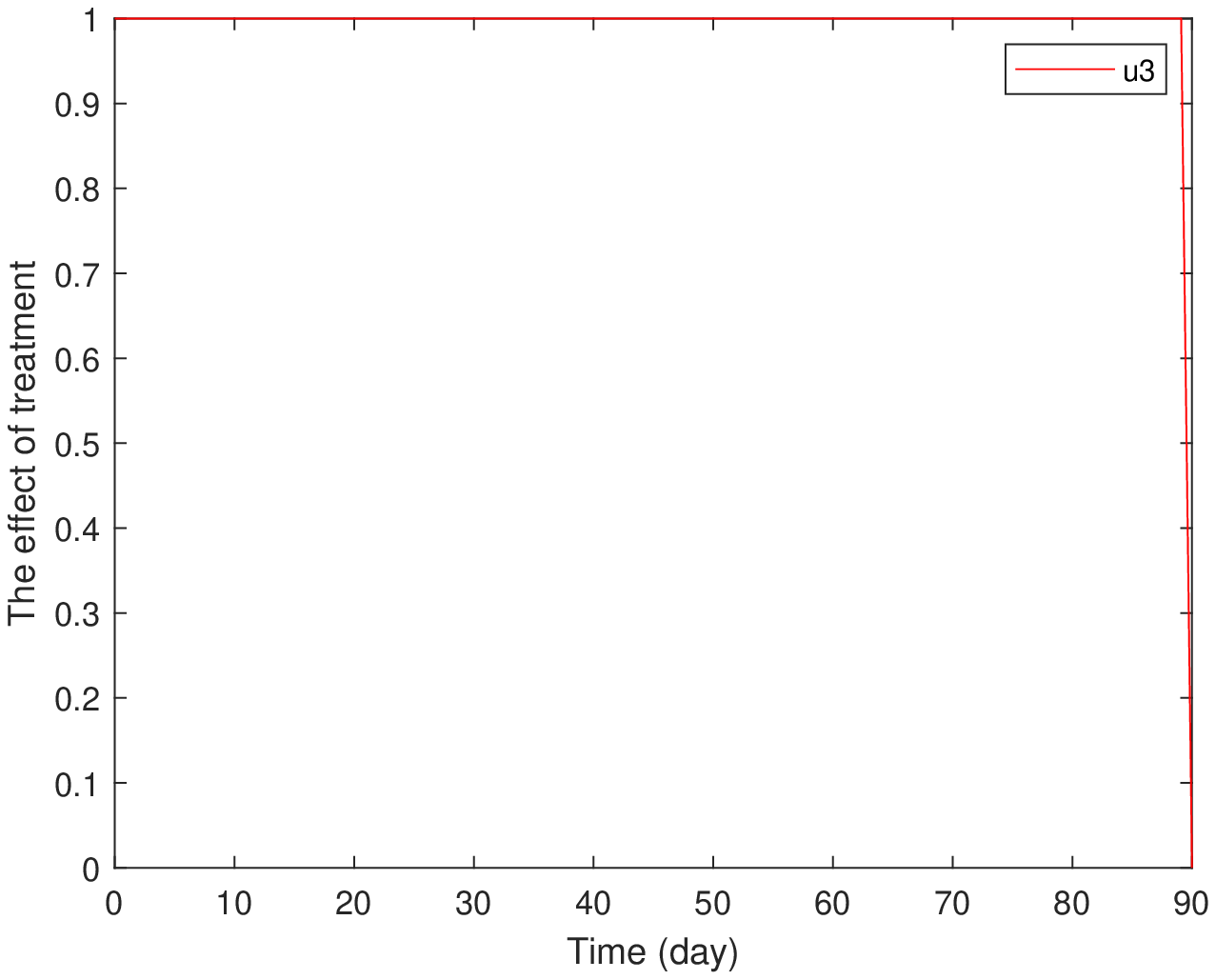}
\caption{$u_3(t)$ (treatment)}
\label{fig2:c3}
\end{subfigure}
\caption{The Pontryagin extremal controls
of the optimal control problem of Section~\ref{sec:4}
for the case of Italy between Sept. 1 and Nov. 30, 2020.}
\label{fig2}
\end{figure*}

In Figures~\ref{fig1} and \ref{fig1:g}, 
we present plots with the numerical solutions 
to the nonlinear differential equations of the generalized SEIR model 
\eqref{q1}, in red color; to the nonlinear differential equations 
of the generalized SEIR control system \eqref{contsys} under optimal controls, 
in the sense of Section~\ref{sec:4}, in green color; 
and the real data of the quarantined cases, 
the number of recovered individuals, and the number 
of deaths from September 1 to October 31, 2020, in orange.
The computed optimal controls for Italy 
from September 1 to November 30, 2020,
which give rise to the green curves in the plots
of Figures~\ref{fig1} and \ref{fig1:g}, are shown in Figure~\ref{fig2}.
The obtained simulations allow us to predict the results 
of the decisions taken in Italy, as well to give the best 
decisions for Italy, according to our generalized SEIR control system
and optimal control problem.

The orange curves in Figures~\ref{fig1:e} and \ref{fig1:f} and Figure~\ref{fig1:g} 
represent the real data on quarantine, recovered, and death cases in Italy from  
September 1 to October 31, 2020. The red curves simulate what 
happens from the beginning of September to the end of November
following the generalized SEIR model \eqref{q1}, 
when the number of quarantined, recovered, and deaths increase, 
and reach, respectively, two million three hundred eighty-eight thousand 
(2388000), nine hundred six thousand three hundred (906300), 
and forty-five thousand seven hundred (45700) cases. 

The red curves in Figures~\ref{fig1:c}, \ref{fig1:d} and \ref{fig1:a} 
simulate what happens from the beginning of September to the end of 
November, according with the generalized SEIR model, when the number 
of infected, exposed and insusceptible cases reach, respectively, 
one million one hundred forty-six thousand (1146000), 
forty-one million two hundred fifty thousand (41250000) 
and five hundred twenty-eight (528) cases.

The green curves in Figures~\ref{fig1} and \ref{fig1:g} show what happens 
from September 1 to November 30, 2020, under optimal control measures,
when the number of infected (Figure~\ref{fig1:c}) 
and recovered (Figure~\ref{fig1:f}) cases increase and 
reach six hundred and fifty (650) 
and two hundred forty-nine thousand four hundred (249400) cases, 
respectively, while the number of exposed (Figure~\ref{fig1:d}), 
insusceptible (Figure~\ref{fig1:a}), and quarantined 
(Figure~\ref{fig1:e}) cases reach eighteen thousand four hundred ninety
(18490), sixty million one hundred eighty thousand (60180000), 
and one hundred twenty-eight (128) cases, respectively. 
Deaths remain stable during the entire period, 
precisely, thirty-five thousand five hundred (35500) 
cases (Figure~\ref{fig1:g}).

The curves in Figure~\ref{fig2} show the optimal controls that need 
to be implemented in order to reduce the overall burden of COVID-19
in Italy and obtain the best possible situation
given by the green curves in Figures~\ref{fig1} and \ref{fig1:g},
which take into account the cost reduction resulting 
from the controls $u_1$, $u_2$ and $u_3$. 
The effect of social distancing is equal to the minimum value of 
its constraint ($u_1=0.1$), see Figure~\ref{fig2:c1}, 
and this corresponds to the application 
of social distancing among the entire population. 
The effect of preventive measures is equal to the maximum value 
of its constraint until September 19 $(u_2=1)$, see Figure~\ref{fig2:c2}, 
then decreases gradually until it reaches zero ($u_2=0$) 
on November 30, 2020, see Figure~\ref{fig2:c3}. 
The effect of treatment takes the maximum value 
of its constraint until November 29, 2020 ($u_3=1$), then decreases 
to zero on November 30, 2020 $(u_3=0)$, meaning
a decrease in the pressure on the health sector.
Note that by taking preventive measures ($u_1$, $u_2$, $u_3$), 
we limit the spread of COVID-19 and we have better results. 
This means that, with the help of optimal control theory, 
what happened in Italy would have been less dramatic.


\section{Conclusion}
\label{sec:conc}

Recent results have shown how the theory of optimal control
is an important tool to combat {COVID-19} 
in a community: in \cite{MyID:459}
for a controlled sanitary deconfinement in Portugal;
in \cite{Zamir} from a more theoretical point of view;
here for the case of Italy.
We proposed a simple SEIR-type control system, 
showing its effectiveness with respect to real data 
from Italy in the period 
from September 1 to November 30, 2020. 
While the real data (see Appendices~A, B and C)
is consistent with the generalized SEIR model \eqref{q1},
because the goal of this model is to describe well the COVID-19 reality, 
our new SEIR control system \eqref{contsys} simulates what would happen 
if we took into account the values of the three 
control functions, as described in Section~\ref{sec:3}. 
In agreement, the situations obtained with controls are better 
than the situations  obtained without controls. More precisely, 
by considering the proposed controls, we show how 
optimal control theory could have drastically
diminish the burden of COVID-19 in the period under study
while taking into account the resulting cost reduction.
In concrete, if it would have been possible
to implement optimally, in the sense of optimal control theory
and Pontryagin's optimality conditions, the control measures of social
distancing as in Figure~\ref{fig2:c1}, preventive means
as in Figure~\ref{fig2:c2}, and treatment as in 
Figure~\ref{fig2:c3}, then it would have been possible
to decrease significantly the number of deaths
(cf. Figure~\ref{fig1:g} and 
Tables~\ref{table:sept:D}~and~\ref{table:oct:D}, 
which account a decrease of 7.36\%
of deaths in Italy by the of October 2020 under optimal control) 
with much less quarantined individuals (see Figure~\ref{fig1:e} 
and Tables~\ref{table:sept:Q} and \ref{table:oct:Q}, 
which account a decrease of 99.96\% of quarantined individuals 
in Italy by the end of October 2020 under optimal control theory).
Thus, one can say that the approach proposed by the theory of optimal control
is very effective, simultaneously from health and economical points of view,
being far from trivial. Note that by following Pontryagin's minimum priciple
one obtains an increase on the number of recovered individuals in a first period,
up to 14-Oct-2020, and, after this date, a decrease on the number of recovered
(cf. Figure~\ref{fig1:f} and Tables~\ref{table:sept:R} and \ref{table:oct:R}),
caused by the drastic reduction on the number of susceptible and infected
(see Figures~\ref{fig1:b} and \ref{fig1:c}, respectively).
While our aim here was to study the effect of controls, 
guided by application of the Pontryagin minimum principle and
showing how they can help to decrease the spread of COVID-19,
other aspects remain open for further research. In particular,
it remains open the theoretical study of the stability of the models.
In this direction, the recent results of 
\cite{MyID:471,F:D:Axioms:2021} may be useful.


\titleformat{\section}{\bf}{\thesection}{1pt}{}

\section*{Acknowledgments}

This research is part of first author's Ph.D. project.
Zaitri is grateful to the financial support from the 
Ministry of Higher Education and Scientific Research of Algeria;
Torres acknowledges the financial support from
CIDMA through project UIDB/04106/2020.
The authors would like to thank two anonymous 
Reviewers for their detailed and thought-out suggestions.



\appendix

\section*{Appendix A: Recovered}

In Tables~\ref{table:sept:R} and \ref{table:oct:R},
we show the real data $R(t)$ of recovered individuals from COVID-19
in Italy, September (Table~\ref{table:sept:R})
and October (Table~\ref{table:oct:R}) 2020, versus the number
$R(t)$ of recovered individuals predicted by SEIR-type model \eqref{q1}
of \cite{Peng} and the controlled model \eqref{contsys}. 
We also indicate the improvement one could have done
by introducing suitable controls, as explained in Section~\ref{sec:3},
and using the theory of optimal control as in Section~\ref{sec:4}.
For that, we give in Tables~\ref{table:sept:R} and \ref{table:oct:R}
the percentage of relative error $\eta_R$ between real data and 
the one predicted by model \eqref{q1}; 
and the improvement $\mathcal{I}_R$ (increase of recovered individuals
with respect to real data in September and up to 14-Oct-2020;
and decrease of recovered from 15-Oct-2020 on, because of a drastic 
reduction on the number of infected and susceptible individuals)  
by introducing controls $u_1$, $u_2$ and $u_3$, as in \eqref{contsys},
in an optimal control way.
\begin{table}	
\small
\caption{Recovered individuals $R(t)$, Sept.~2020.}
\label{table:sept:R}
\begin{center} 	
\begin{tabular}{c c c c c c} \hline
Day & Real & \eqref{q1} & \eqref{contsys} & $\eta_R$ & $\mathcal{I}_R$\\ \hline
01 & 207944 & 207944 & 207944 & 0\% & 0\% \\ 
05 & 209610 & 207996 & 236134 & 0.77\%  & 12.65\% \\
10 & 211885 & 209176 & 238769 & 1.27\%  & 12.68\% \\
15 & 214645 & 211897 & 240170 & 1.28\%  & 11.89\% \\
20 & 218351 & 214873 & 241363 & 1.59\%  & 10.53\% \\
25 & 222716 & 218150 & 242306 & 2.05\%  & 8.79\% \\
30 & 227704 & 221973 & 243000 & 2.51\%  & 6.71\% \\ \hline
\end{tabular}
\end{center}
\end{table}
\begin{table}
\small
\caption{Recovered individuals $R(t)$, Oct.~2020.}
\label{table:oct:R}
\begin{center} 
\begin{tabular}{c c c c c c} \hline
Day & Real & \eqref{q1} & \eqref{contsys} & $\eta_R$ & $\mathcal{I}_R$\\ \hline
01 & 222832 & 224334 & 243132 &0.67\% &9.11\% \\ 
05 & 232681 &  226703     & 243647 &2.56\% &4.71\% \\
10 & 238525 &   232871   & 244263 &2.37\% &2.40\% \\
15 & 245964 &    241255            & 244857 &1.91\% &0.45\% \\
20 & 255005 &    252990       & 245433 &0.79\% &3.75\% \\
25 & 266203 &   269718  & 245994 &1.32\% &7.59\% \\
29 & 279282 &  288247   & 246401 &3.21\% &11.77\% \\ \hline
\end{tabular}
\end{center}
\end{table}


\section*{Appendix B: Deaths}

In Tables~\ref{table:sept:D} and \ref{table:oct:D},
we give real data of death individuals $D(t)$ from COVID-19
in Italy, September (Table~\ref{table:sept:D})
and October (Table~\ref{table:oct:D}) 2020, versus the number
$D(t)$ of death individuals predicted by the SEIR-type model 
\eqref{q1} of \cite{Peng} and our controlled 
model \eqref{contsys}. We also indicate the
improvement one could have done
by introducing suitable controls, as explained in Section~\ref{sec:3},
and using the theory of optimal control as in Section~\ref{sec:4}:
we show the percentage of relative error $\eta_D$ between real data and 
the one predicted by model \eqref{q1}; 
and the improvement $\mathcal{I}_D$ (decrease of death individuals
with respect to real data) by introducing controls 
$u_1$, $u_2$ and $u_3$, as in \eqref{contsys},
in an optimal control way.
\begin{table}
\small	
\caption{Death individuals $D(t)$, Sept.~2020.}
\label{table:sept:D}
\begin{center}
\begin{tabular}{c c c c c c} \hline
Day & Real & \eqref{q1} & \eqref{contsys} & $\eta_D$ & $\mathcal{I}_D$\\ \hline
01 & 35491  & 35491   &  35491 &0\% &0\% \\ 
05 & 35541 &  35510 &  35495 & 0.08\%&0.12\% \\
10 & 35597 &  35538  & 35496 &0.16\% &0.28\% \\
15 & 35645&  35570 & 35496 &0.21\% &0.41\% \\
20 & 35724 &    35606 & 35496 &0.33\% &0.63\% \\
25 & 35818 & 35648    &  35496  &0.47\% &0.89\% \\
30 & 35918 &   35702  &  35497    &0.60\% &1.17\% \\ \hline
\end{tabular}
\end{center}
\end{table}
\begin{table}
\small		
\caption{Death individuals $D(t)$, Oct.~2020.}
\label{table:oct:D}
\begin{center} 
\begin{tabular}{c c c c c c} \hline
Day & Real & \eqref{q1} & \eqref{contsys} & $\eta_D$ & $\mathcal{I}_D$\\ \hline
01 & 35941  & 35715 & 35497 &0\% &0\% \\ 
05 & 36030  &  35773     & 35497 &0.71\%&1.47\% \\
10 & 36166 &   35870   & 35497 &0.81\%&1.84\%\\
15 & 36427 &    36008          & 35497 &1.15\%&2.55\% \\
20 &  36832 &    36206       & 35497 &1.69\%& 3.62\% \\
25 & 37479 &   36491   & 35497 &2.63\%  &5.28\%  \\
29 & 38321 &  37003  &  35498 &3.43\%&7.36\% \\ \hline
\end{tabular}
\end{center}
\end{table}


\section*{Appendix C: Quarantined}

In Tables~\ref{table:sept:Q} and \ref{table:oct:Q},
we show the real data $Q(t)$ of quarantined individuals from COVID-19
in Italy, September (Table~\ref{table:sept:Q})
and October (Table~\ref{table:oct:Q}) 2020, versus the number $Q(t)$
of quarantined individuals predicted by the SEIR-type model \eqref{q1}
of \cite{Peng}  and the one predicted by our 
model \eqref{contsys}. We also indicate 
the improvement one could have done
by introducing suitable controls, 
as explained in Section~\ref{sec:3},
and using the theory of optimal control:
we give the percentage of relative error $\eta_Q$ between real data 
and the one predicted by model \eqref{q1}; 
and the improvement $\mathcal{I}_Q$ (decrease of quarantined individuals
with respect to real data) by introducing controls 
$u_1$, $u_2$ and $u_3$ in an optimal control way.
\begin{table}
\small	
\caption{Quarantined individuals $Q(t)$, Sept.~2020.}
\label{table:sept:Q}
\begin{center}
\begin{tabular}{c c c c c c} \hline
Day & Real & \eqref{q1} & \eqref{contsys} & $\eta_Q$ & $\mathcal{I}_Q$\\ \hline
01  & 26754 & 26754   &  26754 &0 \% &0 \% \\ 
05  & 31194  &  29264 &  1023 &06.18\% &96.69\% \\
10  & 35708   &  31105 & 337 &12.89\%&99.05\%\\
15  & 39712   &  32183 & 228 &18.95\% &99.42\% \\
20  & 44098     &   34808 & 176&21.06\%&99.60\%\\
25  & 47718     &  39848      &  149  &16.49\%&99.68\%\\
30  & 51263     &  48428      &  134   &05.53\%&99.73\%\\ \hline
\end{tabular}
\end{center}
\end{table}
\begin{table}
\small	
\caption{Quarantined individuals $Q(t)$, Oct.~2020.}
\label{table:oct:Q}
\begin{center}
\begin{tabular}{c c c c c c} \hline
Day & Real & \eqref{q1} & \eqref{contsys} & $\eta_Q$ & $\mathcal{I}_Q$\\ \hline
01  & 52647 & 50023 & 130 &04.98\%&99.75\% \\ 
05  & 58903  &  62193    & 124 &05.58\% &99.78\% \\
10  & 74829 &   83557   & 119 &11.63\% &99.84\% \\
15  & 99266 &    116035         & 116 &16.89\% &99.88\% \\
20  &  142739 &   164668      & 112 &15.36\%&99.92\% \\
25  & 222241 &  236520   & 109 &06.42\% &99.95\%\\
29  & 299191&  317055  &  107 &05.97\%&99.96\%\\ \hline
\end{tabular}
\end{center}
\end{table}


\end{document}